\documentclass[journal]{dcls}
\usepackage[utf8]{inputenc}
\usepackage{textcomp}
\usepackage{setspace}
\usepackage[inkscapelatex=false]{svg}

\usepackage{graphicx}
\usepackage{amsmath}
\usepackage{mathtools}
\usepackage[version=4]{mhchem}
\usepackage{siunitx}
\usepackage{longtable,tabularx}
\usepackage{placeins}
\usepackage{makecell}
\usepackage{algorithm}
\usepackage{algpseudocode}
\algrenewcommand\algorithmicrequire{\textbf{Inputs:}}
\algrenewcommand\algorithmicensure{\textbf{Outputs:}}

\usepackage{booktabs}
\usepackage{subcaption}
\usepackage{array}
\usepackage{multirow}
\title{Fuel-optimal Boost-Back Guidance via Successive Convexification with a Terminal Constraint for Reusable Launch Vehicles}

\author[1]{Mingyeom Kim}
\author[1]{Changhyun Kang}
\author[1]{Dohyun Lee}
\author[1,*]{Byeong-Un Jo}

\affil[1]{Department of Aerospace System Engineering, Sejong University, Seoul, 05006, Republic of Korea}

\newcommand{\corrfootnote}{%
\begingroup
\renewcommand{\thefootnote}{*}
\footnotetext{Corresponding author, E-mail: bjo@sejong.ac.kr}
\endgroup
}

\begin{document}
\setstretch{1.2}
\maketitle
\corrfootnote

\begin{abstract}
This paper proposes a fuel-optimal guidance algorithm for the boost-back phase of reusable launch vehicles. We formulate the guidance problem as a free-final-time optimal control problem and solve it by successive convexification (SCvx). The key feature of the proposed formulation is to impose a terminal constraint on the instantaneous impact point (IIP) at the end of the boost-back burn. Since this constraint depends only on the burnout position and velocity, it confines the optimization horizon to the powered phase while enforcing that the predicted ballistic impact point coincides with the target under a spherical-Earth, central-gravity model. When the entire trajectory from stage separation to the landing site is discretized as a single phase, the long ballistic coast must also be propagated through the discretized dynamics, and the resulting defect accumulates over the coast. The bang-off switching point falls between discretization nodes as well. The proposed formulation, in contrast, reduces the problem size and represents the bang-off structure explicitly, with the powered phase resolved by the discretization and the coast handled analytically. The terminal constraint is derived in two forms: 1) the non-iterative, closed-form Keplerian IIP and 2) the F\&G solution, which propagates the orbit with the Lagrange coefficients $F$ and $G$ expressed in terms of the eccentric anomaly. We adopt complex-step differentiation to evaluate the Jacobians of both constraints, for which closed-form expressions are not readily available. We validate the proposed algorithm through case studies on the Falcon 9 CRS-10 mission and its return-to-launch-site (RTLS) scenario. We also examine the trade-off between computational cost and optimality by comparing the proposed algorithm with a single-phase formulation of the entire trajectory, a closed-form guidance law that considers the flight path angle rate, and a benchmark solution obtained by offline trajectory optimization. The results indicate that imposing the IIP constraint analytically provides near-optimal performance with reduced computational cost.
\end{abstract}



\section{Introduction}
Reusable launch vehicles (RLVs) have transformed modern space transportation by enabling the recovery and reuse of the first stage, which drastically reduces launch costs and increases launch frequency\cite{jo2022optimal}. To further reduce launch costs, the recovery maneuvers need to minimize propellant consumption, thereby maximizing the payload capacity. Immediately after stage separation, the first stage performs a boost-back burn that moves the instantaneous impact point (IIP) toward the launch site. Among the three descent burns, namely the boost-back, reentry, and landing burns, the boost-back burn requires the largest total impulse. It therefore has a significant impact on payload capacity and is particularly important in return-to-launch-site (RTLS) missions.\\
Early studies focused on the feasibility of booster return itself, either by jointly optimizing booster sizing and the return trajectory\cite{mckinney1986vehicle} or by comparing trajectory options for returning to the launch site without a secondary propulsion system\cite{hellman2009return}. Pitch-over maneuvers and the associated rocket-back guidance were subsequently studied, and the return trajectory began to be considered as a guidance problem\cite{lu2011pitch,su2011integrated}. Ahn and Roh\cite{ahn2012noniterative} proposed an algorithm that predicts the Keplerian IIP without iteration, and later expressed its time derivative in terms of the disturbing acceleration\cite{ahn2014analytic}. With this relation, the IIP could be actively controlled rather than merely predicted. Jo and Ahn\cite{jo2018near} decomposed the IIP rate into downrange and crossrange components and demonstrated near-optimal performance. However, the performance of this method degrades in RTLS scenarios, in which the IIP must be moved opposite to the downrange direction. To address this, Jo et al. \cite{jo2025fuel} proposed two algorithms: an IIP guidance law that also controls the flight path angle, and a minimum-impulse guidance law based on the hit equation. Both algorithms reduce the onboard computational burden by relying on analytic relations. The former, however, requires the target flight path angle rate to be specified in advance, which limits operational flexibility and leaves a relatively large optimality gap. The latter involves iterative computation. Moreover, neither algorithm can explicitly enforce path constraints such as the throttle range or the attitude limits arising from sensor field-of-view requirements.\\
Such constraints can be handled explicitly by posing the return trajectory directly as an optimization problem. Direct methods that discretize the trajectory using pseudospectral collocation and solve the resulting nonlinear program\cite{patterson2014gpops} can impose path and terminal constraints in their original form. They are therefore widely used to obtain benchmark solutions for trajectory design. For the atmospheric reentry of winged RLVs, Mishra and Sushnigdha\cite{mishra2024novel} constructed a reference trajectory from piecewise polynomials inside a reentry corridor in the height--velocity plane, which is defined by path constraints such as the heat flux and the dynamic pressure. The polynomial parameters that minimize the terminal range-to-go error were then determined through a metaheuristic search. For such nonlinear programming and heuristic methods, however, both convergence and the required computational effort depend on the initial guess and on the problem instance. Neither can be guaranteed a priori, which makes it difficult to meet onboard real-time requirements.\\
Convex optimization, in contrast, has two properties that address this limitation. First, a convex problem solved by an interior-point method converges to the global optimum regardless of the initial guess. Second, an upper bound on the number of interior-point iterations required to reach a prescribed accuracy is known a priori as a polynomial function of the problem size. The convergence behavior is thus inherent in the problem structure rather than dependent on algorithm tuning or on the initial guess. In their pioneering work, A\c{c}{\i}kme\c{s}e and Ploen\cite{acikmese2007convex} formulated powered descent guidance for Mars landing as a second-order cone programming (SOCP) problem. Using Pontryagin's maximum principle, they showed that the nonconvex lower bound on the thrust magnitude can be relaxed without loss of optimality. Sagliano\cite{Sagliano2018PSCP} transcribed this SOCP formulation using pseudospectral differentiation and quadrature operators in place of the conventional trapezoidal scheme, which substantially reduced the error between the optimal states and those obtained by propagating the optimal control history. This approach was later extended to a generalized $hp$ pseudospectral formulation that trades off accuracy against computational time\cite{Sagliano2019hpPSCP}. Despite the favorable properties of convex optimization, the need to re-solve the optimization problem at every guidance cycle had long hindered its onboard use.  Lu\cite{lu2017introducing} observed that the practice of guidance and control is shifting from closed-form equations toward algorithms that rely on intensive onboard computation, and termed this trend computational guidance and control. Since then, real-time solvers customized to the problem structure\cite{dueri2017customized,elango2022customized,kamath2023customized}, together with improvements in onboard computing performance, have made onboard execution practically feasible. Kwon et al. \cite{kwon2021sequential}, for example, demonstrated real-time execution of Mars powered descent guidance on a next-generation onboard processor of the European Space Agency. However, applying convex optimization requires the problem to be expressed in convex form, and the nonlinear dynamics of a launch vehicle preclude a direct SOCP formulation. To handle this nonlinearity, sequential convex programming (SCP) iteratively solves convex subproblems linearized about a reference trajectory. Among SCP methods, the successive convexification (SCvx) algorithm handles both nonlinear dynamics and nonconvex state constraints, and convergence of its iterates to a stationary point has been established for particular trust-region update rules\cite{mao2016successive,mao2018successive}. Szmuk and A\c{c}{\i}kme\c{s}e\cite{szmuk2018successive} extended this approach to the six-degree-of-freedom powered landing problem with nonlinear dynamics and nonconvex constraints. More recent work has further improved the computational efficiency of powered descent guidance. Wang et al. \cite{wang2026bspline}, for example, approximated the state variables with B-splines to reduce the size of the convex subproblems and recovered the control profile analytically after optimization. These studies mainly address the landing burn at relatively low altitudes, whereas the boost-back phase has received little attention. Li et al.\cite{li2021convex} applied SCvx to boost-back trajectory optimization, but they imposed the terminal condition as a target orbit derived from a predetermined atmospheric entry point. The target orbit must therefore be specified in advance, and its right ascension of the ascending node must be updated during the iterations.\\
This paper addresses a boost-back guidance problem in which the terminal target is the IIP near the launch site rather than an atmospheric entry point. The vehicle is modeled with three degrees of freedom in the Earth-centered, Earth-fixed (ECEF) frame. The nonconvex thrust magnitude constraint is handled by lossless convexification\cite{acikmese2007convex} and the nonlinear dynamics by SCvx, with the aim of obtaining the bang-off thrust profile that characterizes the fuel-optimal solution of this class of problems. When the entire trajectory from stage separation to the landing site is treated as a single phase, discretization nodes must also be allocated to the coast phase, which occupies most of the flight time. This increases the problem size and makes the terminal accuracy depend on how accurately the long ballistic coast is propagated by the discretized dynamics. To overcome this, we propose a terminal constraint that confines the optimization horizon to the powered phase and imposes the IIP condition on the state at the end of the boost-back burn, hereafter referred to as burnout. Because this constraint depends only on the burnout position and velocity, the coast phase is excluded from the discretization while the predicted impact point is still driven to the target. The terminal constraint is constructed with two different IIP prediction methods, each used in a separate problem so that their performance can be compared: 1) the noniterative, closed-form solution of Ahn and Roh\cite{ahn2012noniterative}, and 2) the F\&G solution based on the eccentric anomaly. Because the analytic Jacobians of both constraints are cumbersome to derive, the constraints are linearized using complex-step differentiation. A case study of an RTLS scenario based on the Falcon 9 CRS-10 mission shows that the proposed algorithm achieves near-optimal performance, with burn time and propellant consumption nearly identical to those of the benchmark solution obtained by offline trajectory optimization. Compared with the single-phase discretization of the entire trajectory, it reduces both the miss distance and the computational time and represents the bang-off structure explicitly. Although the analytic guidance laws of Jo et al.\cite{jo2025fuel} compute commands rapidly, they exhibit a non-negligible optimality gap in the RTLS scenario. The proposed algorithm solves the fuel-optimal problem directly and thereby reduces this loss of optimality. It also allows thrust and path constraints to be included in the formulation.\\
The contributions of this paper are twofold. First, we propose a formulation that imposes the IIP at burnout as a terminal constraint, which confines the optimization horizon to the powered phase while still driving the predicted impact point to the target. Compared with the single-phase discretization of the entire trajectory, this formulation reduces the number of nodes and the computational time. It also represents the burn-off structure explicitly, with no residual thrust in the coast phase. Second, we derive the terminal constraint in two forms, namely the Keplerian IIP and the eccentric-anomaly-based F\&G solution. These constraints are linearized using complex-step differentiation, which is free of the subtractive cancellation error inherent in finite differences. The proposed formulation also offers several advantages over the preceding closed-form guidance law\cite{jo2025fuel}. The closed-form approach requires little computation because it evaluates an analytic acceleration command at every guidance cycle, but it cannot explicitly impose inequality constraints such as the throttle range or pointing constraints. The proposed formulation, in contrast, can include such constraints directly in the convex subproblems. Furthermore, because the terminal constraint depends only on the burnout position and velocity, it can be replaced by other terminal conditions without changing the structure of the formulation. This property also allows the reentry and landing phases that follow the boost-back burn to be appended as additional phases, so that the entire RTLS trajectory can be optimized as a multi-phase convex optimization problem.\\
The remainder of this paper is organized as follows. Section~II defines the boost-back guidance problem and presents the equations of motion. Section~III describes the linearization, discretization, and SCvx algorithm used to solve the problem. Section~IV derives the two forms of the terminal IIP constraint and their linearization by complex-step differentiation. Section~V presents the RTLS case study, and Section~VI concludes the paper.

\section{Problem Definition}
\subsection{Dynamics}
In this paper, a spherical Earth model is assumed, and the motion of the vehicle is described in the ECEF frame. The three-degree-of-freedom equations of motion of the vehicle in the ECEF frame are
\begin{equation}
\label{eq:dynamics}
\begin{aligned}
\dot{\mathbf{r}}(t) &= \mathbf{v}(t) \\
\dot{\mathbf{v}}(t) &= -\frac{\mu}{{r}(t)^3}\mathbf{r}(t)-\boldsymbol{\Omega}\times(\boldsymbol{\Omega}\times\mathbf{r}(t))-2\boldsymbol{\Omega}\times\mathbf{v}(t)+\frac{\mathbf{T}(t)}{m(t)}\\
\dot{m}(t) &= -\frac{\lVert\mathbf{T}(t)\rVert}{g_0 \cdot I_{\mathrm{sp}}}
\end{aligned}
\end{equation}
where $\mathbf{r}$ is the position vector from the center of the Earth to the vehicle, $\mathbf{v}$ is the velocity vector, $m$ is the mass, and $\mathbf{T}$ is the thrust vector. Here, $\mu$ is the standard gravitational parameter of the Earth, $g_0$ is the gravitational acceleration at sea level, $I_{\mathrm{sp}}$ is the specific impulse, and $\boldsymbol{\Omega}(\coloneqq[0,\:0,\:\omega_E]^\top)$ is the angular velocity vector of the Earth.\\
The nonconvex thrust magnitude constraint is convexified by lossless convexification (LCvx)\cite{acikmese2007convex} through a change of variables and a slack variable. In addition, all variables are normalized to improve the numerical conditioning of the problem. The position and velocity are scaled by the Earth radius $R_E$ and by $\sqrt{R_E g_0}$, respectively, the Earth rotation rate and time are scaled by $\sqrt{g_0/R_E}$ and $\sqrt{R_E/g_0}$, and accelerations are scaled by $g_0$. The same symbols are used for the normalized quantities in Eq.\eqref{eq:dynamics_LCVX} and thereafter. The resulting convexified and normalized dynamics are
\begin{equation}
\label{eq:dynamics_LCVX}
\begin{aligned}
    \dot{\mathbf{r}}(t) &= \mathbf{v}(t) \\
    \dot{\mathbf{v}}(t) &=
    -\frac{1}{{r}(t)^3}\mathbf{r}(t)
    -\boldsymbol{\Omega}\times(\boldsymbol{\Omega}\times\mathbf{r}(t))
    -2\boldsymbol{\Omega}\times\mathbf{v}(t)
    +\mathbf{p}(t) \\
    \dot{z}(t) &= -\alpha\sigma(t)
\end{aligned}
\end{equation}
where $z\coloneqq \ln{m}$, $\mathbf{p}\coloneqq\mathbf{T}/m$, and $\alpha\coloneqq 1/(g_0 I_{\mathrm{sp}})$. Here, $\sigma$ is a slack variable introduced as an upper bound on $\|\mathbf{p}\|$.

\begin{figure}[hbt!]
\centering
\includegraphics[width=0.5\columnwidth]{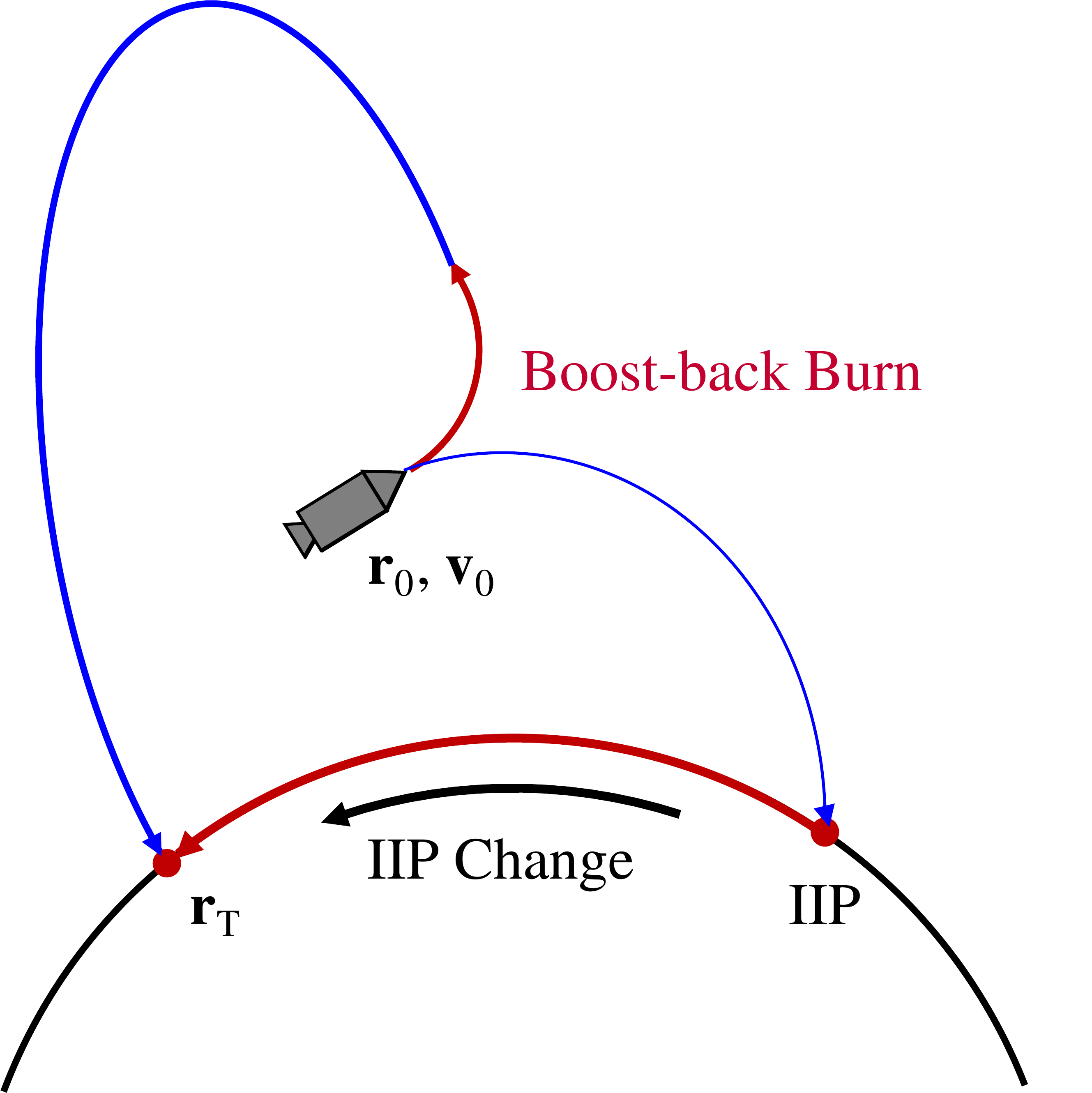}
\caption{Conceptual diagram of the boost-back guidance of an RLV}
\label{Fig_IIP}
\end{figure}

\subsection{Problem Description} \label{subsec:problem1}
Let $\mathbf{r}_0$ and $\mathbf{v}_0$ denote the position and velocity at stage separation. Without thrust, the vehicle follows a free-fall trajectory under gravity alone, as illustrated in Fig.\ref{Fig_IIP}. The IIP is defined as the point at which the vehicle would intersect the Earth's surface under free-fall motion. In the RTLS scenario, the first stage returns to the launch site. In this paper, the goal of the boost-back burn is therefore to place the IIP of the first stage near the launch site. Accordingly, we impose a terminal constraint that requires the IIP at burnout to coincide with the target point $\mathbf{r}_\mathrm{T}$, and we obtain the optimal guidance command through convex optimization. The resulting optimal control problem, Problem~1, is formulated as follows:

\begin{align}
\mathop{\mathrm{minimize}}\limits_{\mathbf{r},\mathbf{v},z,\mathbf{p},\sigma,t_f}\quad
    & -z(t_f) \label{eq:obj}\\[4pt]
\text{subject to} \quad & \nonumber\\[4pt]
\underline{\text{Dynamics:}}\quad
    & \text{Eq.}\eqref{eq:dynamics_LCVX} \nonumber\\
\underline{\text{Path Constraints:}}\quad
    & z_\mathrm{dry} \leq z(t) \leq z_\mathrm{wet}\label{eq:const_z}\\
    & s \geq s_\mathrm{min}\\
    & 0\leq\sigma(t)\leq T_\mathrm{max}e^{-z(t)}\label{eq:const_LCvx1}\\
    & \lVert \mathbf{p}(t)\rVert\leq\sigma(t)\label{eq:const_LCvx2}\\[4pt]
\underline{\text{Terminal Constraint:}}\quad
    & f_\mathrm{IIP}\left(\mathbf{r}(t_f),\mathbf{v}(t_f)\right)=\mathbf{r}_\mathrm{T}
      \label{eq:const_IIP}\\[4pt]
\underline{\text{Boundary Conditions:}}\quad
    & \mathbf{r}(t_0)=\mathbf{r}_0\label{eq:boundary_r}\\
    & \mathbf{v}(t_0)=\mathbf{v}_0\label{eq:boundary_v}\\
    & z(t_0)=z_0\label{eq:boundary_z}
\end{align}
In Problem~1, Eq.\eqref{eq:obj} is the objective function, which maximizes the terminal mass and therefore minimizes the propellant consumption. Equation\eqref{eq:dynamics_LCVX} represents the convexified dynamics, and Eq.\eqref{eq:const_z} bounds the mass within its physically admissible range. Equations\eqref{eq:const_LCvx1}--\eqref{eq:const_LCvx2} are the control constraints convexified by LCvx. Equation\eqref{eq:const_IIP} is the IIP constraint at burnout, which is the key constraint of this paper, and Eqs.\eqref{eq:boundary_r}--\eqref{eq:boundary_z} are the initial conditions.

\section{Successive Convexification for Boost-Back Guidance}
\subsection{Linearization}
To construct a convex subproblem of the nonlinear and nonconvex Problem~1, all of its nonconvex elements must be convexified. With the state $\mathbf{x}(t)=[\mathbf{r}^\top(t),\;\mathbf{v}^\top(t),\;z(t)]^\top$ and the control input $\mathbf{u}(t)=[\mathbf{p}^\top(t),\;\sigma(t)]^\top$, the original nonlinear dynamics are written as
\begin{equation}
\label{eq:xdot(t)=f(x)}
\dot{\mathbf{x}}(t) = f(t,\mathbf{x}(t),\mathbf{u}(t)), \quad t\in [0,t_f]
\end{equation}
The free-final-time problem is converted into a fixed-final-time problem by introducing the normalized time $\tau\in[0,1]$, as shown in Eq.\eqref{eq:Time-dilation}. The dilation factor $s(\coloneqq \mathrm{d}t/\mathrm{d}\tau)$ maps the normalized time back to physical time, so that the free final time is absorbed into a single optimization variable $s$ that multiplies the right-hand side of the dynamics.
\begin{equation}
\label{eq:Time-dilation}
\overset{\circ}{\mathbf{x}}(t) = \frac{\mathrm{d}}{\mathrm{d}\tau}\mathbf{x}(t)=\frac{\mathrm{d}t}{\mathrm{d}\tau}\frac{\mathrm{d}}{\mathrm{d}t}\mathbf{x}(t)=\frac{\mathrm{d}t}{\mathrm{d}\tau}\dot{\mathbf{x}}(t)=s\cdot \dot{\mathbf{x}}(t) \coloneqq F(t,\mathbf{x}(t),\mathbf{u}(t),s)
\end{equation}
A first-order Taylor series expansion about the reference trajectory $(\bar{\mathbf{x}},\bar{\mathbf{u}},\bar{s})$ converts the nonlinear dynamics in Eq.\eqref{eq:Time-dilation} into the linear time-varying (LTV) dynamics in Eq.\eqref{eq:dynamics_LTV}. The coefficient matrices are the Jacobians evaluated along the reference trajectory, as defined in Eq.\eqref{eq:linearization}.
\begin{equation}
\label{eq:dynamics_LTV}
\overset{\circ}{\mathbf{x}}(\tau)=A(\tau)\mathbf{x}(\tau)+B(\tau)\mathbf{u}(\tau)+C(\tau)s+D(\tau)
\end{equation}
where
\begin{subequations}\label{eq:linearization}
\begin{align}
A(\tau) &\coloneqq \nabla_{\mathbf{x}} F(\tau,\bar{\mathbf{x}}(\tau),\bar{\mathbf{u}}(\tau),\bar{s}),\\
B(\tau) &\coloneqq \nabla_{\mathbf{u}} F(\tau,\bar{\mathbf{x}}(\tau),\bar{\mathbf{u}}(\tau),\bar{s}),\\
C(\tau) &\coloneqq \nabla_{s} F(\tau,\bar{\mathbf{x}}(\tau),\bar{\mathbf{u}}(\tau),\bar{s}),\\
D(\tau) &\coloneqq F(\tau,\bar{\mathbf{x}}(\tau),\bar{\mathbf{u}}(\tau),\bar{s}) -A(\tau)\bar{\mathbf{x}}(\tau) -B(\tau)\bar{\mathbf{u}}(\tau) -C(\tau)\bar{s}.
\end{align}
\end{subequations}

\subsection{Discretization}
\begin{figure}[hbt!]
\centering
\includegraphics[width=0.7\columnwidth]{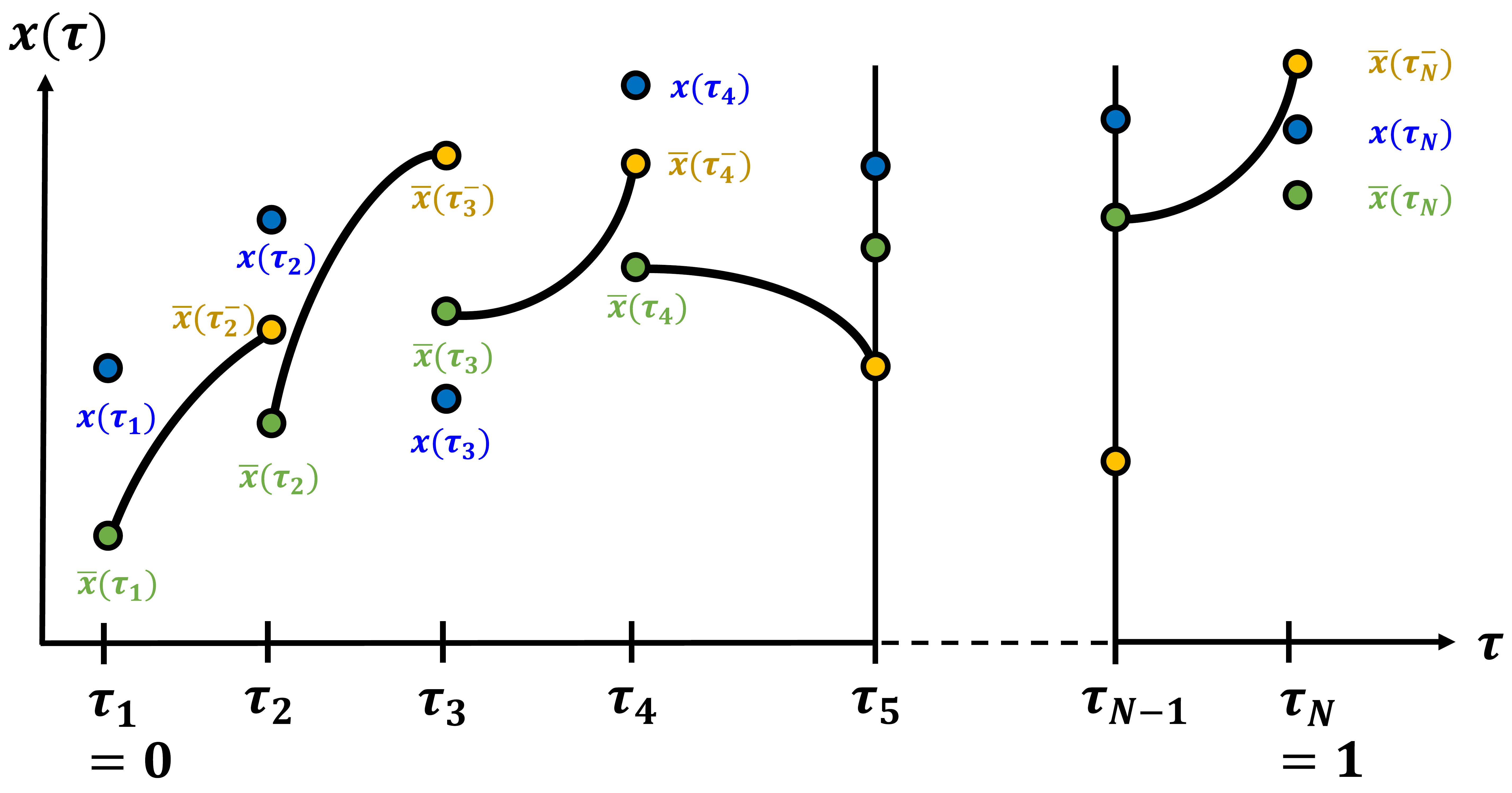}
\caption{Propagation of the state}
\label{Propagation}
\end{figure}

Discretization methods include those that parameterize only the control input, such as zero-order hold and first-order hold (FOH), and pseudospectral methods that parameterize both the state and the control input. Among these, we adopt FOH to discretize Eq.\eqref{eq:dynamics_LTV}, because FOH is known to yield sparsity patterns that significantly reduce computational time
\cite{malyuta2019discretization}. With FOH, the control input is defined only at the temporal nodes and is linearly interpolated between successive nodes, as in Eq.\eqref{u(tau)}.
\begin{equation}\label{u(tau)}
    \mathbf{u}(\tau)=\lambda_{k}^{-}(\tau)\,\mathbf{u}_{k}+\lambda_{k}^{+}(\tau)\,\mathbf{u}_{k+1},\quad \forall\tau\in[\tau_{k},\tau_{k+1}]
\end{equation}
where
\begin{equation*}\label{u_foh}
    \lambda_{k}^{-}(\tau) \coloneqq \frac{\tau_{k+1}-\tau}{\tau_{k+1}-\tau_{k}},\quad\lambda_{k}^{+}(\tau) \coloneqq \frac{\tau-\tau_{k}}{\tau_{k+1}-\tau_{k}},\quad k=1:N-1
\end{equation*}
With the piecewise-affine control in Eq.\eqref{u(tau)}, the LTV dynamics can be expressed in terms of deviations from the reference for all $\tau\in[\tau_k,\,\tau_{k+1}]$, as shown in Eq.\eqref{foh_LTV}\cite{kamath2023customized}. Here, $\Delta\square$ denotes the deviation of a variable from its reference value, i.e., $\Delta\square\coloneqq\square-\bar{\square}$, and $\Delta\overset{\circ}{\mathbf{x}}(\tau)=\overset{\circ}{\mathbf{x}}(\tau)-F(\tau,\bar{\mathbf{x}}(\tau),\bar{\mathbf{u}}(\tau),\bar{s})$.
\begin{equation}\label{foh_LTV}
    \Delta\overset{\circ}{\mathbf{x}}(\tau)=A(\tau)\Delta \mathbf{x}(\tau)+B(\tau)\lambda_{k}^{-}(\tau)\,\Delta \mathbf{u}_{k}+B(\tau)\lambda_{k}^{+}(\tau)\,\Delta \mathbf{u}_{k+1}+C(\tau)\Delta s
\end{equation}
Equation\eqref{foh_LTV} has a unique solution, which is expressed using the state transition matrix (STM) $\Phi(\tau,\tau_k)$ as
\begin{equation}
    \Delta{\mathbf{x}}(\tau)=\Phi(\tau,\tau_k)\Delta \mathbf{x}(\tau_k)+\int_{\tau_k}^{\tau}\Phi(\tau,\zeta)\{B(\zeta)\lambda_k^-(\zeta)\Delta \mathbf{u}_k+B(\zeta)\lambda_k^+(\zeta)\Delta \mathbf{u}_{k+1}+C(\zeta)\Delta s\}d\zeta
\end{equation}
Defining $\Psi_{A}(\tau)\coloneqq\Phi(\tau,\tau_k)$ and applying the \textit{stitching condition}, we recover the discretized dynamics in terms of the absolute variables, as shown in Eq.\eqref{eq:Discretized_Dynamics}.
\begin{equation}\label{eq:Discretized_Dynamics}
     \mathbf{x}_{k+1}=A_{k}\mathbf{x}_k+B_{k}^{-}\mathbf{u}_k+B_{k}^{+}\mathbf{u}_{k+1}+C_{k}s+D_{k}
\end{equation}
where
\begin{equation}
    D_{k} \coloneqq \bar{\mathbf{x}}(\bar{\tau}_{k+1})-(A_{k}\bar{\mathbf{x}}_k+B_{k}^{-}\bar{\mathbf{u}}_k+B_{k}^{+}\bar{\mathbf{u}}_{k+1}+C_{k}\bar{s})
\end{equation}
Because Eq.\eqref{eq:Psi_A} is integrated together with Eq.\eqref{eq:Ak,Bk} over each interval, the FOH coefficient matrices $A_k$, $B_k^-$, $B_k^+$, and $C_k$ can be obtained without matrix inversions.
\begin{subequations}\label{eq:Ak,Bk}
\begin{align}
A_{k} &= I+ \lim_{y\rightarrow\tau_{k+1}^-}\int_{\tau_k}^{y}A(\zeta)\Psi_{A}(\zeta)\,d\zeta\\
B_{k}^{-} &= \lim_{y\rightarrow\tau_{k+1}^-}\int_{\tau_{k}}^{y}\{A(\zeta)\Psi_{B^{-}}(\zeta)+B(\zeta)\lambda_{k}^{-}(\zeta)\}\,d\zeta\\
B_{k}^{+} &= \lim_{y\rightarrow\tau_{k+1}^-}\int_{\tau_k}^{y}\{A(\zeta)\Psi_{B^{+}}(\zeta)+B(\zeta)\lambda_{k}^{+}(\zeta)\}\,d\zeta\\
C_{k} &= \lim_{y\rightarrow\tau_{k+1}^-}\int_{\tau_k}^{y}\{A(\zeta)\Psi_{C}(\zeta)+C(\zeta)\}\,d\zeta
\end{align}
\end{subequations}
where
\begin{subequations}\label{eq:Psi_A}
   \begin{align}
\overset{\circ}{\Psi}_{A}(\tau) &= A(\tau)\Psi_{A}(\tau)\\
\overset{\circ}{\Psi}_{B^{-}}(\tau) &= A(\tau)\Psi_{B^{-}}(\tau)+B(\tau)\lambda_{k}^{-}(\tau)\\
\overset{\circ}{\Psi}_{B^{+}}(\tau) &= A(\tau)\Psi_{B^{+}}(\tau)+B(\tau)\lambda_{k}^{+}(\tau)\\
\overset{\circ}{\Psi}_{C}(\tau) &= A(\tau)\Psi_{C}(\tau)+C(\tau)
   \end{align}
\end{subequations}
These auxiliary matrices are integrated over each interval from the initial conditions
$\Psi_{A}(\tau_k)=I$ and $\Psi_{B^{-}}(\tau_k)=\Psi_{B^{+}}(\tau_k)=\Psi_{C}(\tau_k)=0$,
and $\bar{\mathbf{x}}(\bar{\tau}_{k+1})$ is obtained by propagating Eq.\eqref{eq:Time-dilation} from $\bar{\mathbf{x}}_k$ with the reference control and dilation factor of the previous iteration, as illustrated in Fig.\ref{Propagation}:
\begin{equation}
    \bar{\mathbf{x}}(\bar{\tau}_{k+1})=\bar{\mathbf{x}}_k+\lim_{y\rightarrow\tau_{k+1}^{-}}\int_{\tau_k}^{y} F(\zeta,\bar{\mathbf{x}}(\zeta),\bar{\mathbf{u}}(\zeta),\bar{s})d\zeta
\end{equation}

\subsection{Successive Convexification}
Although the preceding procedures convert the infinite-dimensional nonconvex problem into a finite-dimensional convex subproblem, solving the resulting subproblems naively can lead to two issues. The first is artificial infeasibility, in which linearization renders a subproblem infeasible even though the original problem is feasible. The second is artificial unboundedness, in which linearization renders the cost of a subproblem unbounded from below even though the cost of the original problem is bounded\cite{mao2018successive}.\\
Artificial infeasibility is addressed by adding slack variables to the linearized constraints. A slack variable added to the linearized dynamics is referred to as a virtual control, whereas one added to the other linearized constraints is referred to as a virtual buffer. The IIP always lies on the Earth's surface, so that $\|f_\mathrm{IIP}\|\equiv R_E$ and the Jacobians of $f_\mathrm{IIP}$ have no component in the radial direction. Consequently, the three linearized equations obtained from Eq.\eqref{eq:const_IIP} cannot be satisfied in the radial direction unless the reference IIP coincides with the target point. To remove this artificial infeasibility, we add a virtual buffer $\boldsymbol{\nu}_\mathrm{IIP}\in\mathbb{R}^3$ to the linearized IIP terminal constraint:
\begin{equation}\label{eq:const_lin_IIP}
    f_\mathrm{IIP}(\bar{\mathbf{r}}_N,\bar{\mathbf{v}}_N)+\frac{\partial f_\mathrm{IIP}(\bar{\mathbf{r}}_N,\bar{\mathbf{v}}_N)}{\partial \mathbf{r}}(\mathbf{r}_N-\bar{\mathbf{r}}_N)+\frac{\partial f_\mathrm{IIP}(\bar{\mathbf{r}}_N,\bar{\mathbf{v}}_N)}{\partial \mathbf{v}}(\mathbf{v}_N-\bar{\mathbf{v}}_N)+\boldsymbol{\nu}_\mathrm{IIP}=\mathbf{r}_{\mathrm{T}}
\end{equation}
The virtual buffer allows Eq.\eqref{eq:const_lin_IIP} to hold for arbitrary $\mathbf{r}_N$ and $\mathbf{v}_N$, thereby removing this source of infeasibility. Note that a nonzero $\boldsymbol{\nu}_\mathrm{IIP}$ in the converged solution would violate the original constraint. A 1-norm penalty with a large weight $w_\mathrm{IIP}$ is therefore added to the objective function to drive $\boldsymbol{\nu}_\mathrm{IIP}$ toward zero:
\begin{equation}
    J_\mathrm{vb}=w_\mathrm{IIP}\lVert \boldsymbol{\nu}_\mathrm{IIP} \rVert_{1}
\end{equation}
At convergence the virtual buffer decreases to a negligible magnitude, so the linearized terminal constraint recovers the original nonlinear condition. The linearized constraints, however, are valid only in the vicinity of the reference trajectory, and the iterates may move too far from it. To prevent this, a trust region is imposed at each iteration to keep the solution close enough to the reference trajectory for the linearization to remain valid and, at the same time, to mitigate the artificial unboundedness described earlier. We adopt a soft trust region, which appears as a squared 2-norm penalty in the objective function:
\begin{equation}
    J_\mathrm{tr}=w_\mathrm{tr}\sum_{k=1}^{N}\lVert \mathbf{x}_k-\bar{\mathbf{x}}_k \rVert_{2}^{2} + w_\mathrm{tr_{s}}\lVert s-\bar{s} \rVert_{2}^{2}
\end{equation}
This penalty keeps the state $\mathbf{x}_k$ and the dilation factor $s$ close to their reference values $\bar{\mathbf{x}}_k$ and $\bar{s}$. Unlike a hard trust region, whose radius must be updated explicitly, the soft trust region penalizes the distance from the reference trajectory, so no explicit radius update is required. If $w_{(\cdot)}$ is too small, the solution may move so far from the reference trajectory that the linearization is no longer valid. If it is too large, a suboptimal solution is likely to be obtained\cite{benedikter2021convex}. The weights must therefore be selected appropriately. Combining the preceding results, we obtain the convex subproblem solved at each SCvx iteration, Problem~2:
\begin{align*}
    \mathop{\mathrm{minimize}}\limits_{\mathbf{x}_k,\mathbf{u}_k,s,\boldsymbol{\nu}_\mathrm{IIP}}\quad & -z_{N}+w_{\mathrm{tr}}\sum_{k=1}^{N}\lVert \mathbf{x}_k-\bar{\mathbf{x}}_k \rVert_{2}^{2} + w_{\mathrm{tr_s}}\lVert s-\bar{s} \rVert_{2}^{2} + w_{\mathrm{IIP}}\lVert \boldsymbol{\nu}_{\mathrm{IIP}} \rVert_{1}\\
    \text{subject to}\quad
    & \mathbf{x}_{1}=\mathbf{x}_{0}\\
    & \mathbf{x}_{k+1}=A_{k}\mathbf{x}_k+B_{k}^{-}\mathbf{u}_k+B_{k}^{+}\mathbf{u}_{k+1}+C_{k}s+D_{k} && k=1:N-1\\
    & \lVert \mathbf{p}_k \rVert_{2}\leq\sigma_{k} && k=1:N\\
    & 0\leq\sigma_{k}\leq T_{\mathrm{max}}e^{-\bar{z}_{k}}(1+\bar{z}_{k}-z_{k}) && k=1:N\\
    & z_{\mathrm{dry}} \leq z_{k} \leq z_{\mathrm{wet}} && k=1:N\\
    & f_{\mathrm{IIP}}(\bar{\mathbf{r}}_N,\bar{\mathbf{v}}_N)+\frac{\partial f_{\mathrm{IIP}}(\bar{\mathbf{r}}_N,\bar{\mathbf{v}}_N)}{\partial \mathbf{r}}(\mathbf{r}_N-\bar{\mathbf{r}}_N)+\frac{\partial f_{\mathrm{IIP}}(\bar{\mathbf{r}}_N,\bar{\mathbf{v}}_N)}{\partial \mathbf{v}}(\mathbf{v}_N-\bar{\mathbf{v}}_N)+\boldsymbol{\nu}_{\mathrm{IIP}}=\mathbf{r}_\mathrm{T}
 \end{align*}
Problem~2 is obtained from Problem~1 by linearizing the nonlinear dynamics and the nonconvex thrust upper bound about the reference trajectory, adding a virtual buffer to the nonlinear IIP constraint, and adding penalty terms to the objective function. Because the discretization nodes $k=1$ and $k=N$ correspond to stage separation and burnout, respectively, the initial conditions of Problem~1 are imposed at the first node. The constraints of Problem~2 consist of affine equalities, affine inequalities, and second-order cones, and the objective augments the linear cost with a quadratic trust-region penalty and a 1-norm virtual-buffer penalty. Problem~2 is therefore a second-order cone program. Its solution is used as the reference for the next iteration, and the process is repeated until the relative change of the objective value between successive iterations falls below $\varepsilon_\mathrm{tol}$. The linearization of $f_{\mathrm{IIP}}$ is described in Section\ref{sec:terminal_constraints}.

\section{Terminal Constraints}
\label{sec:terminal_constraints}
If Problem~1 in Section\ref{subsec:problem1} is solved with Eq.\eqref{eq:const_IIP} replaced by $\mathbf{r}(t_f)=\mathbf{r}_\mathrm{T}$, treating the trajectory from stage separation to the landing site as a single phase and allocating nodes over the entire horizon, the optimal solution should exhibit a bang-off structure. However, when the LCvx-relaxed problem is discretized and solved numerically, small residual thrust remains in the coast phase. More importantly, discretizing the coast requires the ballistic trajectory to be integrated node by node, and the resulting defect accumulates over the entire coast phase. In this paper, the optimization horizon is therefore confined to the powered phase of the boost-back burn, and the coast phase after burnout is accounted for through the terminal constraint in Eq.\eqref{eq:const_lin_IIP}. This approach removes the coast phase from the optimization problem, so that the ballistic trajectory is propagated analytically rather than node by node, while concentrating the $N$ nodes on the powered phase to reduce both the problem size and the computational time.\\
Equation\eqref{eq:const_lin_IIP} is based on three assumptions. First, only central gravity acts on the vehicle after burnout. Second, the Earth is a sphere. Third, the IIP is the first intersection of the trajectory with the Earth's surface. Section\ref{subsec:keplerian_IIP} derives the Keplerian IIP constraint from the noniterative, closed-form solution of Ahn and Roh\cite{ahn2012noniterative}. Section\ref{subsec:FGsol} then derives the F\&G solution constraint, which propagates the trajectory to the IIP with the Lagrange coefficients $F$ and $G$ by solving for the eccentric anomaly at impact. In Sections\ref{subsec:keplerian_IIP} and \ref{subsec:FGsol}, the subscript 0 denotes the burnout state and should be distinguished from the initial conditions in Table\ref{tab:parameter}.

\subsection{Keplerian IIP}\label{subsec:keplerian_IIP}
The impact point is computed in an inertial frame that is instantaneously aligned with the ECEF frame at burnout. The burnout velocity is converted by $\mathbf{v}_I=\mathbf{v}_E+\boldsymbol{\Omega}\times\mathbf{r}_E$, the Keplerian impact point is obtained in that frame, and the result is rotated back through $\omega_E t_p$. Because the frame is defined at burnout, only the coast duration $t_p$ enters the rotation. Let $\mathbf{i}_{\mathbf{r}0}$ and $\mathbf{i}{\mathbf{v}_0}$ denote the unit vectors of the burnout position and velocity in this frame. The IIP is then expressed in terms of the flight path angle $\gamma_0$ and the flight angle $\phi$ as in Eq.\eqref{eq:hit}, where $\phi$ is the angle traveled by the vehicle from burnout to the IIP.
\begin{equation}\label{eq:hit}
     \mathbf{i}_p=\frac{\mathbf{r}_p}{r_p}=\frac{\cos(\phi+\gamma_0)}{\cos\gamma_0}\mathbf{i}_{\mathbf{r}_0}+\frac{\sin\phi}{\cos\gamma_0}\mathbf{i}_{\mathbf{v}_0}
\end{equation}
The flight angle $\phi$ is given by
\begin{equation}
    \phi=\arcsin\left(\frac{c_1c_3+\sqrt{c_1^2c_3^2-(c_1^2+c_2^2)(c_3^2-c_2^2)}}{(c_1^2+c_2^2)}\right),
\end{equation}
where
\begin{equation}
    c_1=-\tan\gamma_0,\quad c_2=1-\frac{1}{\lambda\cos^2\gamma_0},\quad c_3=\frac{r_0}{r_p}-\frac{1}{\lambda\cos^2\gamma_0}
\end{equation}
Equation\eqref{eq:tp_kepler} gives the time of flight from burnout to the IIP, which is required to correct the IIP for the Earth's rotation over the flight.
\begin{equation}\label{eq:tp_kepler}
    t_p=\frac{r_0}{v_0\cos\gamma_0}\left(\frac{\tan\gamma_0(1-\cos\phi)+(1-\lambda)\sin\phi}{(2-\lambda)\left(\frac{1-\cos\phi}{\lambda\cos^2\gamma_0}+\frac{\cos(\gamma_0+\phi)}{\cos\gamma_0}\right)}+\frac{2\cos\gamma_0}{\lambda\left(\frac{2}{\lambda}-1\right)^{1.5}}\arctan\left(\frac{\sqrt{\frac{2}{\lambda}-1}}{\cos\gamma_0\cot\frac{\phi}{2}-\sin\gamma_0}\right)\right),
\end{equation}
where
\begin{equation*}
    \lambda \coloneqq \frac{r_0v_0^2}{\mu}
\end{equation*}

\subsection{F\&G Solution}\label{subsec:FGsol}
Whereas the previous subsection computes the IIP directly in closed form in terms of $\phi$, this subsection propagates the coast trajectory in terms of the change in eccentric anomaly $\Delta E$. This change is determined from the condition that the trajectory intersects the Earth's surface. Here, $\mathbf{r}_0$ and $\mathbf{v}_0$ denote the burnout position and velocity vectors in the ECI frame. Equation\eqref{Eccen_rp} gives the radius after a change $\Delta E$ from burnout, where the semi-major axis $a$ is obtained from the vis-viva equation and $\sigma_0$ from the dot product of the position and velocity vectors. Setting $r_p(\Delta E)=R_E$ and solving Eq.\eqref{Eccen_rp} for $\Delta E$ with the Newton--Raphson method yields the eccentric anomaly from burnout to impact.
\begin{equation}\label{Eccen_rp}
    r_p(\Delta{E})=a+(r_0-a)\cos\Delta{E}+\sqrt{a}\sigma_0\sin\Delta{E}
\end{equation}
Substituting $\Delta E$ into Eq.\eqref{eq:tp_fg} yields the time of flight of the ballistic coast. As with Eq.\eqref{eq:tp_kepler}, this value determines the rotation angle required to transform the IIP from the ECI frame to the ECEF frame.
\begin{equation}\label{eq:tp_fg}
    t_{p}=\sqrt{\frac{a^3}{\mu}}\left[\Delta{E}-\frac{\sigma_0}{\sqrt{a}}(\cos\Delta{E}-1)-\left(1-\frac{r_0}{a}\right)\sin\Delta{E}\right],
\end{equation}
where
\begin{equation*}
    a=\frac{1}{\frac{2}{r_0}-\frac{{v}_0^2}{\mu}
},\quad \sigma_0\coloneqq \frac{\mathbf{r}_0 \cdot \mathbf{v}_0}{\sqrt{\mu}}
\end{equation*}
Once $\Delta E$ and $t_{p}$ are determined, the Lagrange coefficients $F$ and $G$ are computed from Eqs.\eqref{eq:F_coeff}--\eqref{eq:G_coeff}, and the position vector of the IIP is expressed as a linear combination of the burnout position and velocity vectors:
\begin{equation}\label{eq:fgsol}
    \mathbf{r}_p=F\mathbf{r}_0+G\mathbf{v}_0,
\end{equation}
where
\begin{align}
    F &= 1-\frac{a}{r_0}(1-\cos\Delta{E}) \label{eq:F_coeff}\\
    G &= t_{p}-\sqrt{\frac{a^3}{\mu}}\left(\Delta{E}-\sin\Delta{E}\right) \label{eq:G_coeff}
\end{align}
Because the methods in Sections\ref{subsec:keplerian_IIP} and \ref{subsec:FGsol} both return the IIP from the burnout position and velocity alone, either can be used directly in the terminal constraint in Eq.\eqref{eq:const_lin_IIP}.

\subsection{Complex-Step Differentiation}\label{subsec:Complex-step}
To include the IIP constraint in the convex subproblem, the constraint must be linearized through a first-order Taylor series expansion about the reference terminal state $(\bar{\mathbf{r}}_N,\bar{\mathbf{v}}_N)$. This expansion requires the Jacobians $\frac{\partial f_{\mathrm{IIP}}(\bar{\mathbf{r}}_N,\bar{\mathbf{v}}_N)}{\partial \mathbf{r}}$ and $\frac{\partial f_{\mathrm{IIP}}(\bar{\mathbf{r}}_N,\bar{\mathbf{v}}_N)}{\partial \mathbf{v}}$ in Eq.\eqref{eq:const_lin_IIP}. Although the Keplerian IIP is available in closed form, its expression is lengthy, and the F\&G solution defines $\Delta E$ implicitly through the Newton--Raphson iteration. In both cases, deriving the analytic Jacobians is cumbersome and error-prone. Finite differences, on the other hand, are limited in accuracy by truncation and round-off errors. We therefore compute the Jacobians using complex-step differentiation, which evaluates the function $f$ with a small imaginary perturbation $ih$ in the complex domain. A Taylor series expansion of the perturbed function yields
\begin{equation*}
    f(x+ih)\simeq f(x)+ihf'(x)-\frac{h^2}{2!}f''(x)-\frac{ih^3}{3!}f'''(x)+\cdots
\end{equation*}
Taking the imaginary part of both sides gives an approximation of the first derivative:
\begin{equation*}
    f'(x)=\frac{\mathrm{Im}[f(x+ih)]}{h}+\mathcal{O}(h^2)
\end{equation*}
The two $3\times3$ Jacobian matrices evaluated at the reference terminal state are
\begin{align}
    J_{\mathbf{r}} &\coloneqq \left.\frac{\partial f_\mathrm{IIP}}{\partial \mathbf{r}}\right|_{\bar{\mathbf{r}}_N,\bar{\mathbf{v}}_N}\in\mathbb{R}^{3\times3}\\
    J_{\mathbf{v}} &\coloneqq \left.\frac{\partial f_\mathrm{IIP}}{\partial \mathbf{v}}\right|_{\bar{\mathbf{r}}_N,\bar{\mathbf{v}}_N}\in\mathbb{R}^{3\times3}
\end{align}
Complex-step differentiation constructs each column of these matrices independently by applying a purely imaginary perturbation to each of the three components of the position and velocity in turn:
\begin{align}
    J_{\mathbf{r}}{(:,j)}&=\frac{\mathrm{Im}(f_{\mathrm{IIP}}(\bar{\mathbf{r}}_N+ih\mathbf{e}_j,\bar{\mathbf{v}}_N))}{h},\quad j=1,2,3\\
    J_{\mathbf{v}}{(:,j)}&=\frac{\mathrm{Im}(f_{\mathrm{IIP}}(\bar{\mathbf{r}}_N,\bar{\mathbf{v}}_N+ih\mathbf{e}_j))}{h},\quad j=1,2,3
\end{align}
where $\mathbf{e}_j$ is the $j$th standard basis vector, and $f_\mathrm{IIP}$ can be either of the functions in Sections\ref{subsec:keplerian_IIP} and \ref{subsec:FGsol}. Complex-step differentiation requires $f_\mathrm{IIP}$ to be analytic in the perturbed variable. All vector magnitudes are therefore evaluated as $\sqrt{\mathbf{x}^\top\mathbf{x}}$ with a non-conjugate transpose, and the built-in absolute-value, norm, and conjugate operations are avoided\cite{martins2003complexstep}. A step size of $h=10^{-20}$ is used. For the F\&G formulation the Newton--Raphson iteration is started from a fixed real initial value $\Delta E=0.5$ and is carried out entirely in complex arithmetic, with at most 100 iterations, and the termination test uses only the real part of the Newton step so that the imaginary perturbation does not affect the stopping decision.

\begin{algorithm}
\caption{Complex-step Linearization of the Terminal Constraint}
\label{alg:csd}
\begin{algorithmic}[1]
\Require
Reference state $\bar{\mathbf{r}}_N,\bar{\mathbf{v}}_N$, step size $h$
\For{$j=1,\ldots,3$}
    \State $\mathbf r^{+} \gets \bar{\mathbf{r}}_N+i h\,\mathbf e_j$
    \State $\displaystyle\frac{\partial f}{\partial \mathbf{r}_j}\gets\frac{\mathrm{Im}\left[f_{\mathrm{IIP}}
    (\mathbf r^{+}, \bar{\mathbf v}_N)\right]}{h}$
    \State $\mathbf v^{+}\gets\bar{\mathbf{v}}_N+i h\,\mathbf e_j$
    \State $\displaystyle\frac{\partial f}{\partial \mathbf{v}_j}\gets\frac{\mathrm{Im}\left[f_{\mathrm{IIP}}
    (\bar{\mathbf r}_N,\mathbf v^{+})\right]}{h}$
\EndFor
\Ensure
$\displaystyle\frac{\partial f_{\mathrm{IIP}}(\bar{\mathbf{r}}_N,\bar{\mathbf{v}}_N)}{\partial \mathbf{r}},\;\frac{\partial f_{\mathrm{IIP}}(\bar{\mathbf{r}}_N,\bar{\mathbf{v}}_N)}{\partial \mathbf{v}}$
\end{algorithmic}
\end{algorithm}

\section{Case Study}
We consider the boost-back phase of the Falcon 9 first stage in the CRS-10 mission, which performed an RTLS landing. The IIP formulations in Sections\ref{subsec:keplerian_IIP} and \ref{subsec:FGsol} are both derived in the ECI frame, whereas the dynamics in Eq.\eqref{eq:dynamics} are defined in the ECEF frame. The IIP must therefore be corrected for the Earth's rotation over the time of flight $t_p$. The rotation matrix from the ECI frame to the ECEF frame, $R_{E/I}$, is
\[R_{E/I} =
\begin{bmatrix}
\cos(\omega_E\,t_p)  & \sin(\omega_E\,t_p) & 0\\
-\sin(\omega_E\,t_p) & \cos(\omega_E\,t_p) & 0\\
0                    &  0                  & 1
\end{bmatrix}\]

\begin{table}[ht]
\centering
\setlength{\tabcolsep}{12pt}
\caption{Parameters}\label{tab:parameter}
\begin{tabularx}{0.8\linewidth}{XXX}
\toprule
\toprule
\textbf{Parameter} & \textbf{Value} & \textbf{Units}\\
\midrule
$I_{\text{sp}}$                              & $311$                     & s\\
$T_{\text{max}}$                             & $279.6$                   & tonf\\
$\text{Target IIP}\,(\text{Lat},\text{Lon})$ & $[28.486,\;-80.543]$      & deg\\
$\mathbf{r}_0$                               & $[1051,\;-5541,\;3160]$   & km \\
$\mathbf{v}_0$                               & $[1.475,\;0.098,\;1.395]$ & km/s\\
$m_0$                                        & $79.2$                    & ton \\
$m_{\text{dry}}$                             & $22.2$                    & ton \\
$s_\mathrm{min}$                             & $31$                      & s \\
$N$                                          & $43$                      & --\\
$\varepsilon_\mathrm{tol}$                   & $1e-4$                    & --\\
$[w_{\text{tr}},\:w_{\text{tr}_\text{s}},\:w_{\text{IIP}}]$ & [1e-3,\:1e-2,\:1e+2] & --\\
\bottomrule
\bottomrule
\end{tabularx}
\end{table}

\begin{figure}[hbt!]
\centering
\includegraphics[width=0.5\columnwidth]{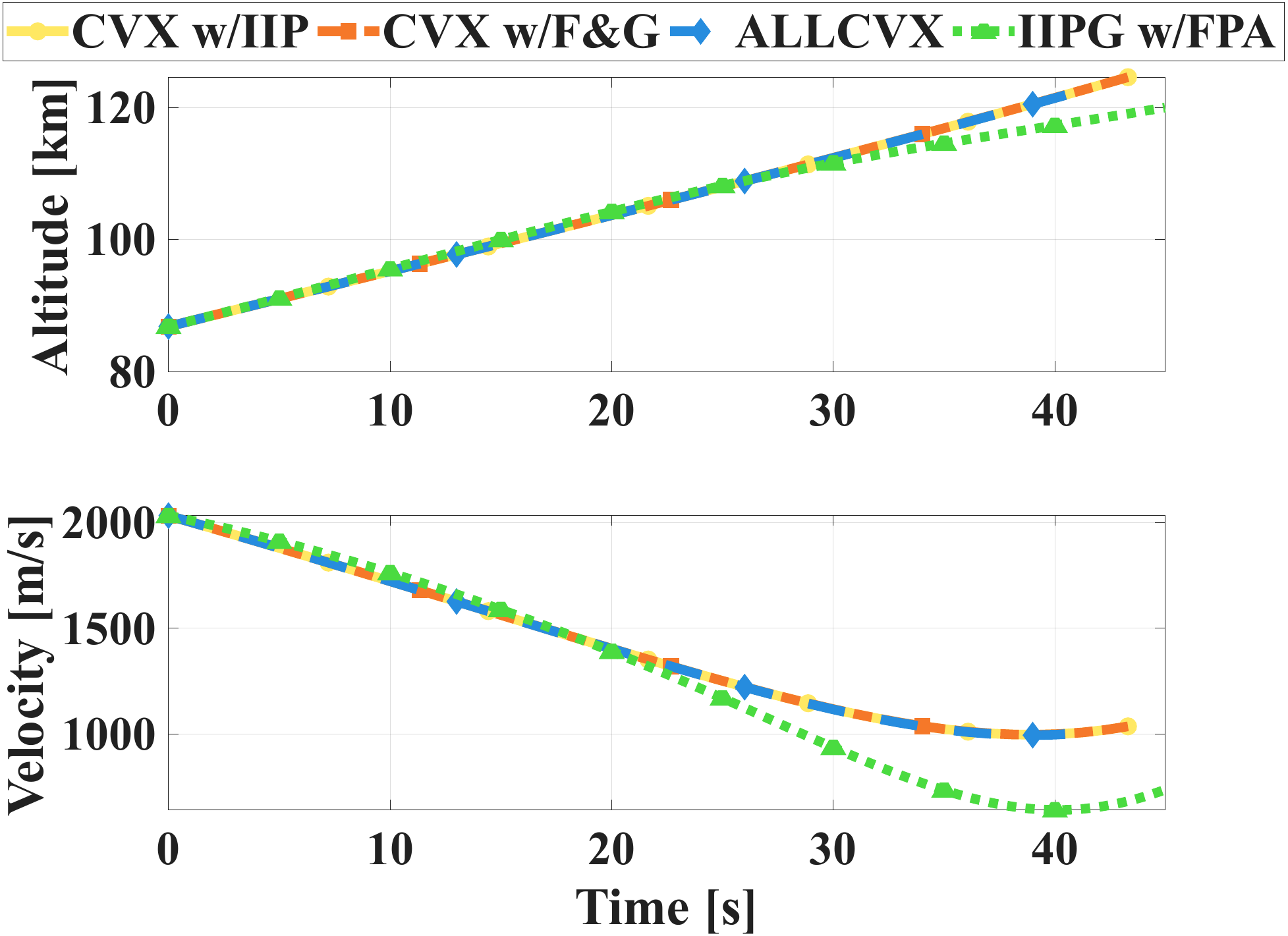}
\caption{Altitude and velocity of the reusable stage in the ECI frame}
\label{Fig:Alt_Vel}
\end{figure}

\begin{figure}[hbt!]
\centering
\includegraphics[width=0.5\columnwidth]{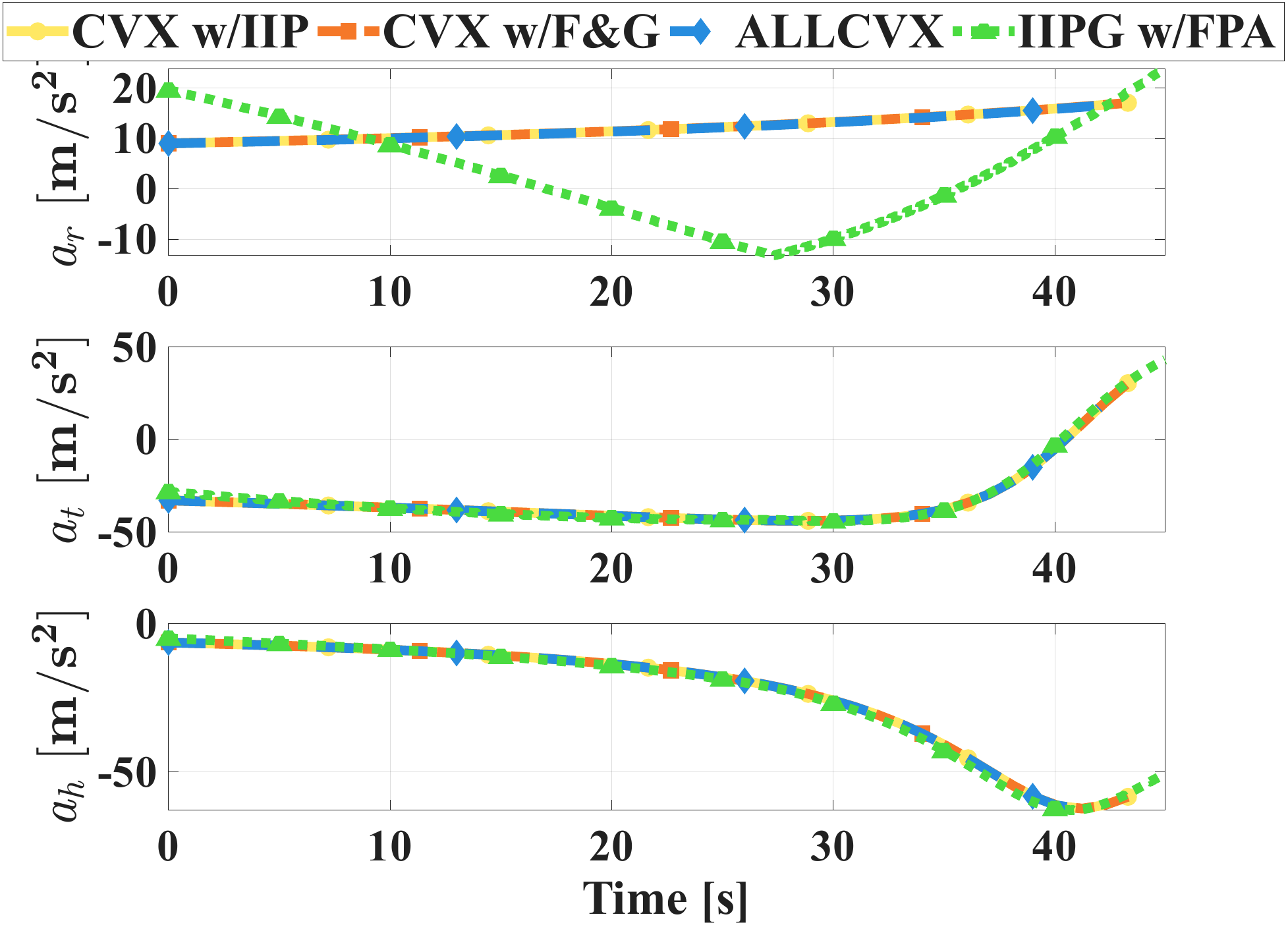}
\caption{Acceleration command profiles}
\label{Fig:Accel}
\end{figure}

\begin{figure}[hbt!]
\centering
\includegraphics[width=0.5\columnwidth]{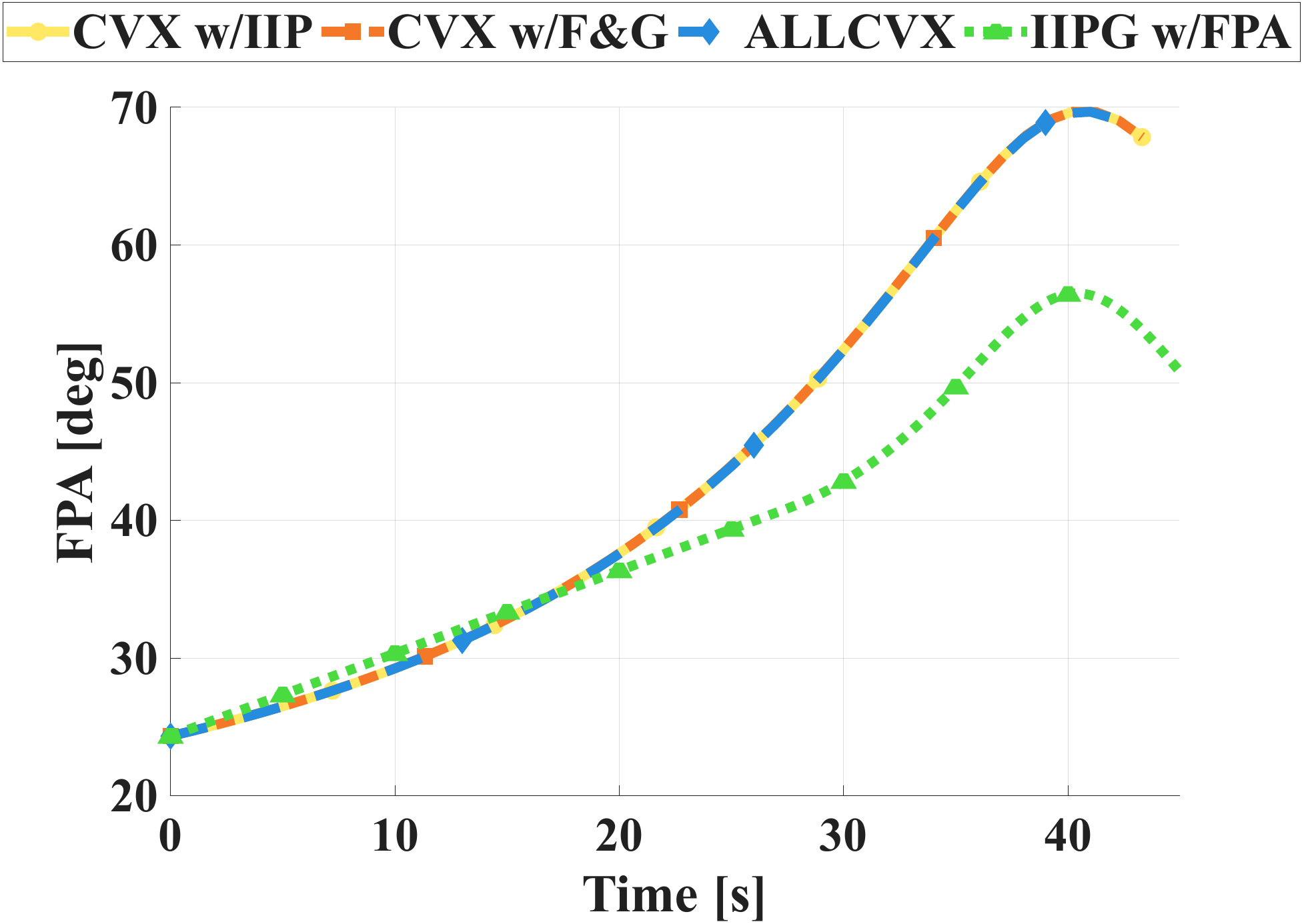}
\caption{Flight path angle of the reusable stage}
\label{Fig:FPA}
\end{figure}

\begin{figure}[hbt!]
\centering
\includegraphics[width=0.5\columnwidth]{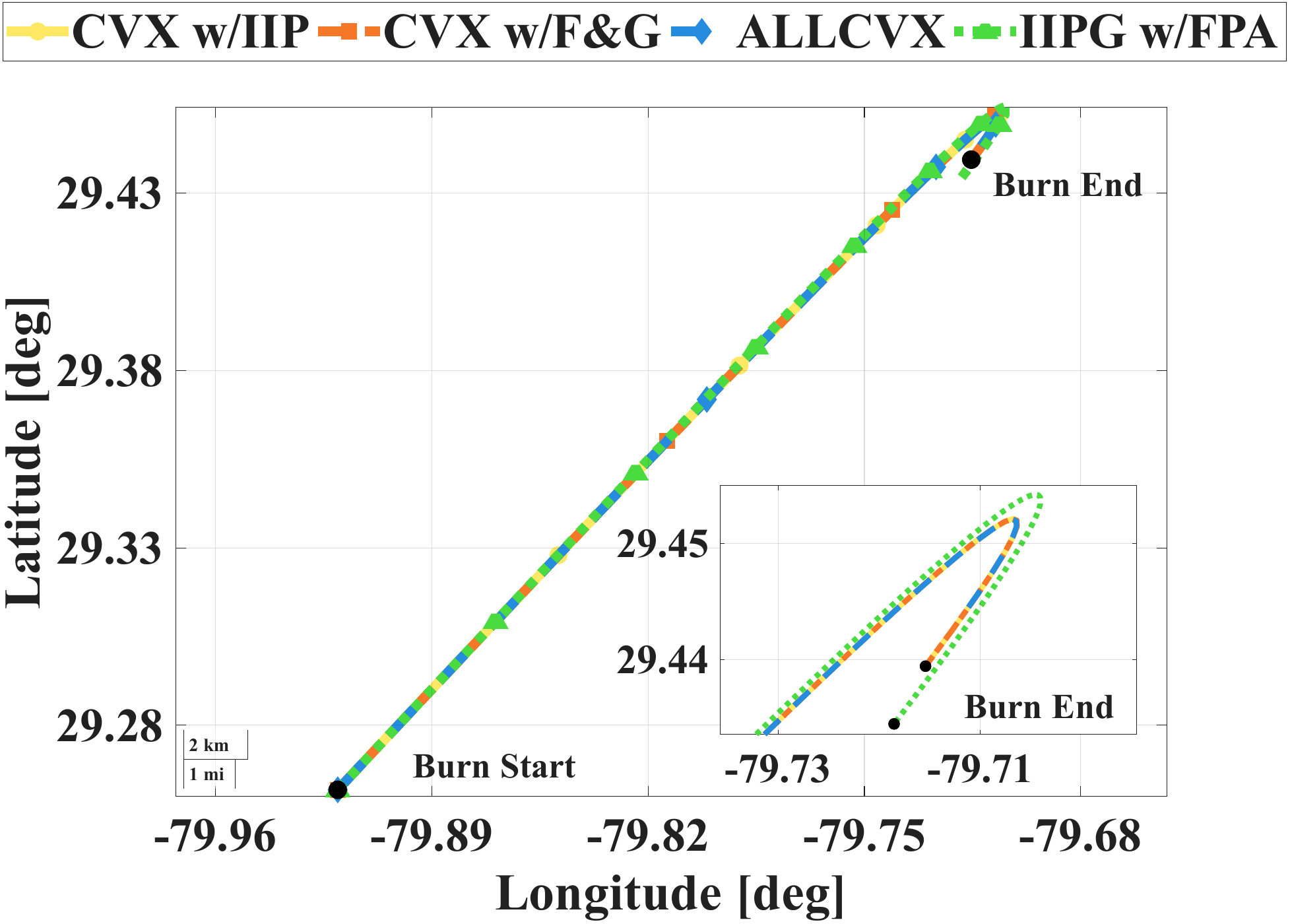}
\caption{Ground track of the reusable stage}
\label{Fig:Geoplot}
\end{figure}

\begin{figure}[hbt!]
\centering
\includegraphics[width=0.5\columnwidth]{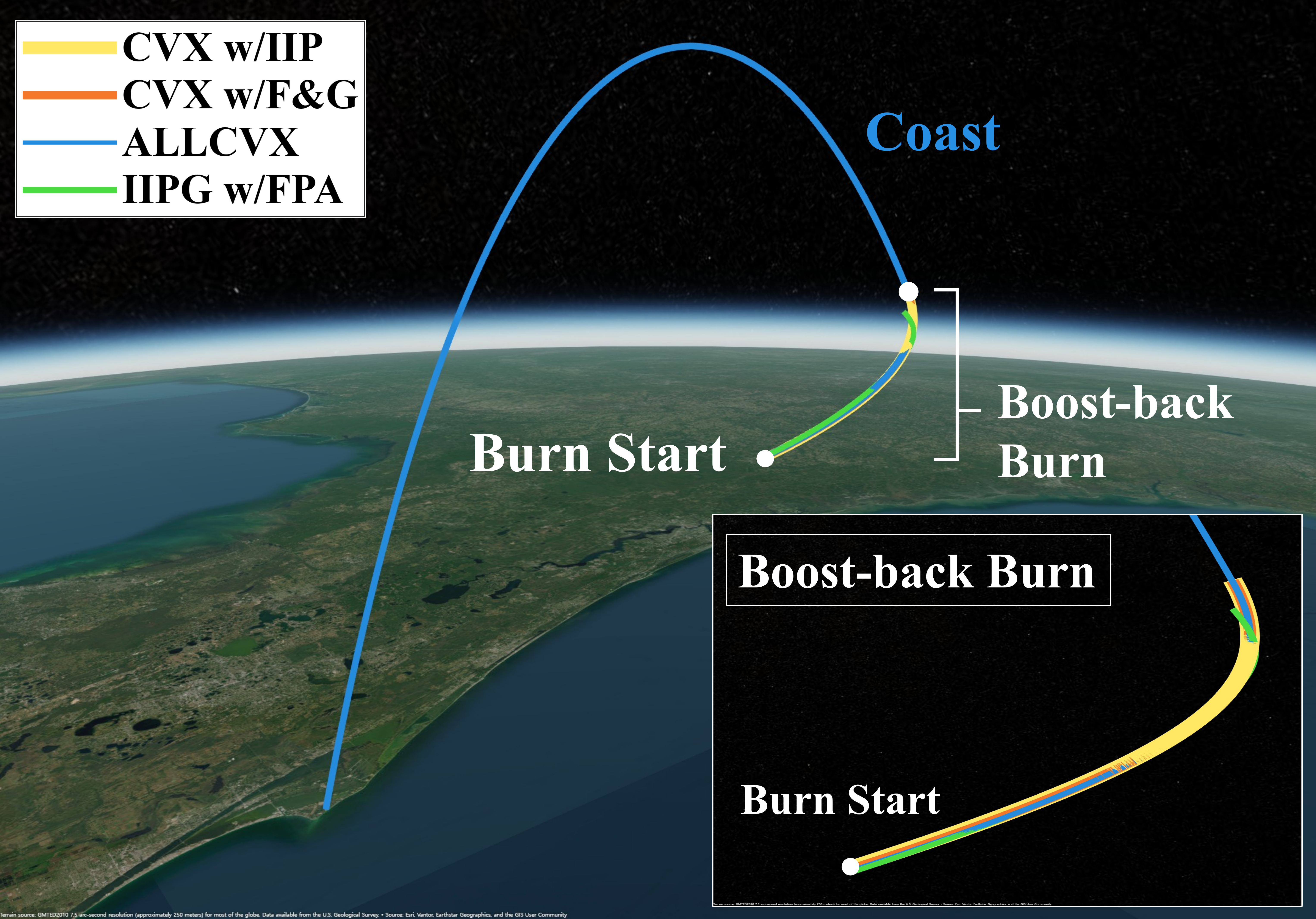}
\caption{Three-dimensional trajectory of the reusable stage}
\label{Fig:3D_Trajectory}
\end{figure}

Table\ref{tab:parameter} lists the vehicle specifications and initial conditions used to solve Problem~2, together with the parameter values and the stopping criterion of the SCvx algorithm. The vehicle specifications and initial conditions are based on publicly available data for the CRS-10 mission, and the initial conditions are given in the ECI frame. Because the ECI and ECEF frames coincide at liftoff and the boost-back burn starts 160~s later, the initial conditions are converted to the ECEF frame by rotating them through the Earth rotation angle over this interval\cite{jo2025fuel}.\\
To examine the effects of different IIP constraints on guidance performance and computational time, we compare four boost-back guidance methods under the same objective function. Three of them share the SCvx-based trajectory optimization framework and differ only in how the landing condition is formulated. The first, CVX w/IIP, computes the Keplerian IIP analytically from the burnout position and velocity and imposes the terminal constraint in Eq.\eqref{eq:const_IIP}. The second, CVX w/F\&G, propagates the coast phase with the Lagrange coefficients $F$ and $G$ and imposes the same constraint at the end of the propagated trajectory. The third, ALLCVX, allocates discretization nodes over the entire trajectory, including the coast phase, and constrains the position at the final node to coincide with the target point. The fourth, IIPG w/FPA, is the noniterative, closed-form IIP guidance law considering the flight path angle proposed by Jo et al.\cite{jo2025fuel} and is referred to hereafter as the conventional method.\\
Figs.\ref{Fig:Alt_Vel}--\ref{Fig:3D_Trajectory} present the results for the RTLS scenario. Figs.\ref{Fig:Alt_Vel}--\ref{Fig:FPA} show the altitude and velocity profiles, acceleration command profiles in the LVLH frame, and flight path angle profiles in the ECI frame during the boost-back burn, respectively. In the LVLH frame, $a_r$, $a_t$, and $a_h$ denote the radial, transverse, and out-of-plane components of the acceleration command. The radial component is positive outward from the center of the Earth, the out-of-plane component is positive along the angular momentum vector, and the transverse component completes the right-handed triad, pointing in the direction of motion. Figs.\ref{Fig:Geoplot} and \ref{Fig:3D_Trajectory} show the ground track during the boost-back burn in the ECEF frame and the trajectory from the boost-back burn through the coast phase to the target point, visualized using the MATLAB Mapping Toolbox.

As shown in Figs.\ref{Fig:Alt_Vel}--\ref{Fig:3D_Trajectory}, the trajectories of the three SCvx-based methods are nearly identical regardless of how the IIP constraint is formulated, whereas the trajectory of the conventional method differs noticeably. The difference can be attributed to the acceleration command profiles shown in Fig.\ref{Fig:Accel}. The SCvx-based methods determine the thrust direction that minimizes propellant consumption over the entire boost-back phase. IIPG w/FPA, in contrast, computes a closed-form command at every guidance cycle that simultaneously satisfies the shortest-path condition and the specified flight path angle rate\cite{jo2025fuel}.

\begin{table}[ht]
\centering
\caption{Comparison of guidance methods}
\label{tab:comparison}
\renewcommand{\arraystretch}{1.2}
\begin{tabular}{c c c c c c}
\toprule
\toprule
\multicolumn{2}{c}{Method}
  & \makecell{Burn Time,\\ s}
  & \makecell{Propellant\\ Consumption, ton}
  & \makecell{Miss Distance,\\ m}
  & \makecell{Computational\\ Time, s} \\
\midrule
\multirow{3}{*}{SCvx (Proposed)}
& Keplerian IIP & 43.29 & 38.92 & 0.0094 & 1.12 \\
& F\&G Solution & 43.29 & 38.92 & 0.0094 & 1.15 \\
& ALLCVX        & 44.00 & 38.93 & 1.1290 & 1.27 \\
\midrule
\multicolumn{2}{c}{Conventional Method}
& 44.96 & 40.39 & -- & -- \\
\multicolumn{2}{c}{Offline Trajectory Optimization}
& 43.29 & 38.92 & -- & -- \\
\bottomrule
\bottomrule
\end{tabular}
\end{table}

Table\ref{tab:comparison} summarizes the results of each guidance method. All algorithms were implemented in MATLAB R2025b using CVX with the MOSEK solver, and all computations were performed on a desktop computer with a 12-core processor (4.4~GHz) and 32~GB of RAM. CVX w/IIP and CVX w/F\&G converge to the same solution to the precision shown. In contrast, ALLCVX yields a slightly longer burn time, a much larger miss distance, and a far larger problem size. Because CVX w/IIP and CVX w/F\&G treat the coast phase analytically, discretizing only the burn phase is sufficient, and $N=43$ nodes are used over a burn time of about 43~s. ALLCVX must also discretize the coast phase and uses $N=340$ nodes over a total flight time of about 339~s. The node counts are chosen so that the time step is approximately 1~s in both cases. The difference in problem size therefore comes only from the length of the discretized horizon, which for ALLCVX is larger by nearly an order of magnitude. Although the objective function of Problem~1 depends only on the terminal mass and would therefore also admit a solution in which the burn is followed by a coast segment, the throttle in the solution obtained here is at its upper bound at every one of the $N=43$ nodes. The final node is therefore the actual burnout point.\\
In addition, a node with an intermediate thrust level appears at the bang-off switching point, because the switching time does not coincide with a discretization node. The burn time of ALLCVX therefore exceeds that of the other two methods by 0.71~s, although its propellant consumption differs by only 0.01~ton. The larger miss distance of ALLCVX arises from a different source. Its terminal node satisfies the target condition to within the solver tolerance, but re-propagating the coast phase with the true dynamics from the first node at which the thrust vanishes reveals the defect accumulated over the 294~s coast, together with the small residual thrust that this propagation discards. The proposed formulation removes this error source entirely, because the coast is propagated analytically by the Keplerian relations instead of being integrated node by node. These results show that treating the IIP constraint analytically reduces both the computational burden and the miss distance while maintaining comparable guidance performance.\\
The conventional method is advantageous in terms of computational efficiency, and ALLCVX has the advantage of not requiring a separate analytic IIP constraint. CVX w/IIP and CVX w/F\&G converge in four iterations and require 1.12~s and 1.15~s in total, whereas ALLCVX converges in three iterations but requires 1.27~s. The cost per iteration is therefore 0.28~s for CVX w/IIP and 0.42~s for ALLCVX, an increase of 51\%, which follows from the larger number of nodes. Because ALLCVX needs one fewer iteration, the advantage in total computational time is smaller, at 11.8\%. CVX w/F\&G is marginally slower than CVX w/IIP because it solves for the eccentric anomaly with the Newton--Raphson method. The optimality of each method is assessed against a benchmark solution obtained by offline trajectory optimization using GPOPS-II\cite{patterson2014gpops}. 

\FloatBarrier
\section{Conclusion}
In this paper, we presented a fuel-optimal boost-back guidance algorithm that imposes a terminal IIP constraint and solves the resulting problem by SCvx. Because the IIP constraint depends only on the burnout state, the optimization horizon was confined to the powered phase while the predicted ballistic impact point was still driven to the target. As a result, the problem size was reduced, and the bang-off structure was represented explicitly, with no residual thrust in the coast phase. Two forms of the terminal constraint, the Keplerian IIP and the F\&G solution, were considered. Because their analytic Jacobians are cumbersome to derive, complex-step differentiation was used to obtain the Jacobians required for linearization.\\
In the RTLS case study based on the Falcon 9 CRS-10 mission, both terminal constraint methods converged in four iterations. They achieved a burn time of 43.29~s, corresponding to an optimality gap of less than 0.01\% relative to the benchmark solution obtained by offline trajectory optimization (43.29~s). In contrast, ALLCVX yielded a burn time of 44.00~s, which is 1.64\% longer. This difference reflects the misalignment between the discretization nodes and the switching time rather than a loss of fuel optimality, because the two solutions differ in propellant consumption by only 0.01ton. Its miss distance of 1.1290~m, however, is approximately 120 times larger than the 0.0094~m of the proposed methods. Although ALLCVX converged in three iterations rather than four, its total computational time was 1.27~s against 1.12~s for CVX w/IIP, and its cost per iteration was 51\% higher. This increase is attributed to the larger problem size caused by allocating discretization nodes to the coast phase. The conventional method yielded a burn time of 44.96~s and a propellant consumption of 40.39~ton, corresponding to gaps of 3.86\% and 3.78\%, respectively, relative to the benchmark solution.\\
These results show that imposing the IIP constraint analytically improves both the miss distance and the computational efficiency while maintaining fuel optimality comparable to that of the benchmark solution. Between the two formulations, the Keplerian IIP constraint provides lower computational cost while achieving the same guidance performance as the F\&G formulation. Several directions remain for future work. These include extending the optimization to the entire RTLS trajectory, formulating the problem with six degrees of freedom based on dual quaternions, using deep reinforcement learning\cite{liu2026deep} to generate the initial reference trajectory, and developing a dedicated solver for onboard real-time implementation.

\section*{Acknowledgments}
This work was supported by the National Research Foundation of Korea (NRF) grant funded by the Korea government (MSIT) (RS-2025-24683783) and the Information Technology Research Center (ITRC) (IITP-2026-RS-00437494).

\bibliography{sample}

\end{document}